\documentclass[11pt, doublespaced, leqno,noinfoline]{article}

\RequirePackage[OT1]{fontenc}

\RequirePackage[aop]{imsart}

\usepackage{amsmath,amsthm,amssymb,graphicx,natbib,color}

\startlocaldefs
\def\@journal{\ }
\newcommand{\eqdef}{\stackrel{d}{=}}
\newcommand{\probconv}{\stackrel{p}{\longrightarrow}}

\newtheorem{thm}{Theorem}[section]
\newtheorem{remark}[thm]{Remark}
\newtheorem{propn}[thm]{Proposition}

\newtheorem{example}{Example}[section]

\newtheorem{lemma}[thm]{Lemma}
\newtheorem{cor}[thm]{Corollary}

\numberwithin{equation}{section}

\endlocaldefs

\begin{document}

\begin{frontmatter}

\title{STATIONARY SYMMETRIC $\alpha$-STABLE DISCRETE PARAMETER RANDOM FIELDS}

\runtitle{Stable Random Fields}


\author{\fnms{Parthanil} \snm{Roy}\ead[label=e1]{pr72@cornell.edu}}
\address{School of Operations Research\\ and Industrial Engineering\\Cornell University\\Ithaca, NY 14853\\\printead{e1}}\and
\author{\fnms{Gennady} \snm{Samorodnitsky}\ead[label=e2]{gennady@orie.cornell.edu}}
\address{School of Operations Research\\ and Industrial Engineering\\Cornell University\\Ithaca, NY 14853\\\printead{e2}}\thanksref{t1}
\thankstext{t1}{Supported in part by NSF
grant DMS-0303493, NSA grant MSPF-05G-049
 and NSF training grant ``Graduate and Postdoctoral
Training in Probability and Its Applications'' at Cornell
University.}

\affiliation{Cornell University}

\runauthor{P.~Roy and G.~Samorodnitsky}

\begin{abstract}

We establish a connection between the structure of a stationary
symmetric $\alpha$-stable random field $(0 < \alpha < 2)$ and
ergodic theory of non-singular group actions, elaborating on a
previous work by \cite{rosinski:2000}. With
the help of this connection, we study the extreme values of the field
over increasing boxes. Depending on the ergodic theoretical and group
theoretical structures of the underlying action, we observe different
kinds of asymptotic behavior of this  sequence of extreme values.

\end{abstract}

\begin{keyword}[class=AMS]
\kwd[Primary ]{60G60}
\kwd[; secondary ]{37A40}
\end{keyword}

\begin{keyword}
\kwd{random field, stable process, ergodic theory, maxima, extreme
 value theory, group action, non-singular map, dissipative, conservative}
\end{keyword}

\end{frontmatter}

\section{Introduction}

In this paper we study the structure of stationary symmetric
$\alpha$-stable discrete parameter non-Gaussian random fields. A
random field $\{X_t\}_{t\in \mathbb{Z}^d}$ is called a symmetric
$\alpha$-stable ($S\alpha S$) random field if for all $c_1, c_2,
\ldots, c_k \in \mathbb{R}$, and, $t_1,t_2,\ldots,t_k \in
\mathbb{Z}^d$, $\sum_{j=1}^k c_j X_{t_j}$ follows a symmetric
$\alpha$-stable distribution. In this paper we will concentrate on
the non-Gaussian case, and hence, we will assume $0 < \alpha < 2$,
unless mentioned otherwise. For further reference on $S \alpha S$
distributions and processes the reader is suggested to read
\cite{samorodnitsky:taqqu:1994}. A random field $\{X_t\}_{t\in
\mathbb{Z}^d}$ is called stationary if
\begin{equation}
\label{e:stat}
\{X_t\} \eqdef \{X_{t+s}\} \;\;\; \text{for all} \  s \in \mathbb{Z}^d\;.\\
\end{equation}
Stationarity means that the law of the random field is invariant
under the action of the group of shift transformations on the
index-parameter $t \in \mathbb{Z}^d$.

More generally, if $(G,+)$ is a countable abelian group with identity
element $0$, then a random field $\{X_t\}_{t\in G}$ is called
$G$-stationary if \eqref{e:stat} holds for all $s\in G$. Most of the
structure results in this paper have immediate analogs for
$G$-stationary fields. We will mention these briefly along the
way. Even though our main interest lies with $\mathbb{Z}^d$-indexed
random fields, at a certain point in the paper a more general group
structure will become important.

Our first task in this paper is to establish a connection between
ergodic theory of nonsingular $\mathbb{Z}^d$-actions (see Section
$1.6$ of \cite{aaronson:1997}) and $S\alpha S$ random fields. Using
the language of the Hopf decomposition of nonsingular flows a
decomposition of stationary $S\alpha S$ processes was established in
\cite{rosinski:1995}. For a general $d>1$ a similar decomposition of
$S\alpha S$ random fields into independent components was given in
\cite{rosinski:2000}. We show the connection between this
decomposition and ergodic theory. This is done in Section
$\ref{sec2}$, using an approach different from the one-dimensional
case, namely, without referring to the Chacon-Ornstein theorem,
which is unavailable in the case $d>1$.

We use the connection with ergodic theory to study the rate of
growth of the partial maxima sequence $\{M_n\}$ of the random field
$X_t$ as $t$ runs over a $d$-dimensional hypercube of size with an
increasing edge length $n$. In the case $d=1$ it has been shown in
\cite{samorodnitsky:2004a} that this rate drops from $n^{1/\alpha}$
to something smaller as the flow generating the process changes from
dissipative to conservative. One can argue that this phase
transition qualifies as a transition between short and long memory.
In this paper we  establish a similar phase transition result for a
general $d\geq 1$.

In Section $\ref{b_n_section}$, we first discuss the asymptotic
behavior of a certain deterministic sequence which controls the size
of the partial maxima sequence $\{M_n\}$. The treatment here is
different from the one-dimensional case due to unavailability of
Maharam extension theorem (see Theorem $2$ in \cite{maharam:1964})
in the case $d>1$. In this section, we also calculate the rate of
growth of partial maxima of the random field. We show that the rate
of growth of $M_n$ is equal to $n^{d/\alpha}$ if the group action
has a nontrivial dissipative component, and is strictly smaller than
that otherwise.

We discuss connections with the group theoretical properties of the
action in Section \ref{sec_group_theory}. For $S\alpha S$ random
fields generated by conservative actions, we view the underlying
action as a group of nonsingular transformations and study the
algebraic structure of this group to get better estimates on the
rate of growth of the partial maxima. Examples illustrating how the
maxima of a random field can grow are discussed in Section
$\ref{examples}$.


\section{Some Ergodic Theory}\label{sec1}

The details on the notions introduced in this section can be found, for
example, in \cite{aaronson:1997}. Unless stated otherwise, the
statements about sets (e.g. equality or disjointness of  two sets) are
understood as holding up to a set of measure zero with respect to the
underlying measure.

Suppose $(S, \mathcal{S}, \mu)$ is a $\sigma$-finite
standard measure space and $(G,+)$ is a countable group with
identity element $0$. A collection of measurable maps $\phi_t:S
\rightarrow S,\;t \in G$ is called a group action of $G$ on $S$ if
\begin{enumerate}
\item $\phi_0$ is the identity map on $S$, and,
\item $\phi_{u+v}=\phi_u \circ \phi_v$  for all $u, v \in G$\,.
\end{enumerate}
A group action $\{\phi_t\}_{t \in G}$ of $G$ on $S$ is called
nonsingular if $ \mu \circ \phi_t \sim \mu$ for all $t \in G\,$.\\

\noindent A set $W \in \mathcal{S}$ is called a wandering set for
the action $\{\phi_t\}_{t \in G}$ if $\{\phi_t(W):\;t\in G\}$ is a
pairwise disjoint collection. The following result (see Proposition
$1.6.1$ of \cite{aaronson:1997}) gives a decomposition of $S$ into two
disjoint and invariant parts.

\begin{propn}
Suppose $G$ is a countable group and $\{\phi_t\}$ is a nonsingular
action of G on $S$. Then $S=\mathcal{C} \cup \mathcal{D}$ where
$\mathcal{C}$
and $\mathcal{D}$ are disjoint and invariant measurable sets  such that\\

1. $\displaystyle{\mathcal{D} = \bigcup_{t \in G} \phi_t(W_{\ast})}$
for some wandering set $W_{\ast}$\,,\\

2. $\mathcal{C}$ has no wandering subset of positive measure.
\end{propn}

\noindent $\mathcal{D}$ is called the dissipative part, and
$\mathcal{C}$  the conservative part of the action. The
action $\{\phi_t\}$ is called
conservative if $S=\mathcal{C}$ and dissipative if $S=\mathcal{D}$\,.\\

\noindent An action $\{\phi_t\}_{t \in G}$ is free if $\mu\big(\{s
\in S:\;\phi_t(s)=s\}\big)=0$ for all $t \in G-\{0\}$. Note that
this definition makes sense because $(S, \mathcal{S})$ is a standard
Borel space  and hence $\{s \in S:\;\phi_t(s)=s\} \in \mathcal{S}$.
The following result is a version of Halmos' Recurrence Theorem for
a nonsingular action of a countable group.

\begin{propn} \label{halmos_rec_thm} Let $\{\phi_t\}$ be a nonsingular
action of a countable group G. If $A \in \mathcal{S}$ and $A \subseteq
\mathcal{C}$, then
\[
\sum_{t \in G} I_A \circ \phi_t = \infty \mbox{ a.e. on }A.
\]

\end{propn}

\begin{proof} Define,
\[
F:=\{s \in S: \exists \, t \in G, t \neq 0 \mbox{ such that
}\phi_t(s)=s\}\,.
\]
Observe that $F$ is $\{\phi_t\}$-invariant. Restrict $\{\phi_t\}$ to
$S-F$. Let $\mathcal{C}_1$ be the conservative part of the
restriction. It is easy to observe that $A \cap F^{c} \subseteq
\mathcal{C}_1$ for all $A \subseteq \mathcal{C}$. Since the
restricted action is free by Proposition $1.6.2$ of
\cite{aaronson:1997} we have,
\begin{eqnarray*}
\sum_{t \in G} I_A \circ \phi_t \geq \sum_{t \in G} I_{A \cap F^{c}}
\circ \phi_t = \infty \mbox{ a.e. on }A \cap F^{c}.
\end{eqnarray*}
Clearly,
\[
\sum_{t \in G} I_A \circ \phi_t = \infty \mbox{ a.e. on }A \cap F.
\]
This completes the proof.
\end{proof}

Recall that the dual operator of a nonsingular transformation $T$ on
$S$ is a linear operator $\hat{T}$ on $L^1(S,\mu)$ such that
\[
\int_S \hat{T}f.g d\mu=\int_S f.g \circ T d\mu\;\;\;\; \text{for all
  $f \in L^1(\mu)$ and $g \in L^\infty(\mu)$\,.}
\]
In particular, if $T$ is invertible, then
\[
\hat{T}f=\frac{d\mu \circ T^{-1}}{d\mu} f \circ T^{-1}\;\;\;\;
\text{for all} \  f \in L^1(\mu) \,,
\]
see Section 1.3 in \cite{aaronson:1997}. The following proposition
is an extension of Theorem $1.6.3$ of \cite{aaronson:1997} to not
necessarily measure-preserving transformations, and can be
established using an argument parallel to that of Propositon 1.3.1
in \cite{aaronson:1997}.

\begin{propn} \label{dualthm}

If $G$ is a countable group and $\{\phi_t\}$ is a nonsingular action
of $G$ on $S$ then for all $f\in L^1(\mu),\, f>0$,
\[
\mathcal{C}=\{s \in S: \sum_{t \in G}\hat{\phi}_t
f(s)=\infty\}\,.
\]
\end{propn}

The following is an immediate corollary, particularly suitable for
our purposes.

\begin{cor} \label{dualcor}
If $G$ is a countable group and $\{\phi_t\}$ is a nonsingular action
of $G$ then
\[
[\sum_{t \in G}\frac{d\mu\circ \phi_t}{d\mu}f\circ
\phi_t=\infty]=\mathcal{C}\;\ \text{for all}\  f\in L^1(\mu),\,
f>0.
\]

\end{cor}

Note that, as mentioned earlier, the equalities of sets in Proposition $\ref{dualthm}$ and Corollary $\ref{dualcor}$ above hold up to sets of $\mu$-measure zero.


\section{Stationary Symmetric Stable Random Fields}\label{sec2}

Suppose $\mathbf{X}=\{X_t\}_{t \in \mathbb{Z}^d}$ is a $S\alpha S$
random field, $0 < \alpha < 2$. We know from Theorem $13.1.2$ of
\cite{samorodnitsky:taqqu:1994} that it has an integral representation of the from
\begin{eqnarray}
X_t&\eqdef& \int_{S} f_t(s)M(ds),\;\; t \in \mathbb{Z}^d\,,
\label{eqn3}
\end{eqnarray}
where $M$ is a $S\alpha S$ random measure on some standard Borel
space $(S,\mathcal{S})$ with $\sigma$-finite control measure $\mu$
and $f_t \in L^\alpha(S,\mu)$ for all $t \in \mathbb{Z}^d$. Note
that $f_t$'s are deterministic functions and hence all the
randomness of $\mathbf{X}$ is hidden in the random measure $M$, and,
the inter-dependence of the $X_t$'s is captured in $\{f_t\}$. The
representation (\ref{eqn3}) is called an integral representation of
$\{X_t\}$. Without loss of generality we can also assume that the
family $\{f_t\}$ satisfies the full support assumption
\begin{eqnarray}
\mbox{Support}\Bigl( f_t,\, t \in
\mathbb{Z}^d\Bigr)=S\,, \label{eqn5}
\end{eqnarray}
because, if that is not the case, we can replace $S$ by
$S_0=\mbox{Support}\bigl( f_t,\, t \in \mathbb{Z}^d\bigr)$ in
$(\ref{eqn3})$.\\

\noindent If, further, $\{X_t\}$ is stationary, then the fact that
the action of the group $\mathbb{Z}^d$ on $\{X_t\}_{t \in
\mathbb{Z}^d}$ by translation of indices preserves the law, and
certian rigidity of spaces $L^\alpha,\, \alpha<2$ guarantees
existence of intergral representations of a special form. This has
been established in \cite{rosinski:1995} for $d=1$ and
\cite{rosinski:2000} for a general $d$. Specifically, there always
exists a representation of the form
\begin{eqnarray}
f_t(s)&=&c_t(s){\left(\frac{d \mu \circ \phi_t}{d
\mu}(s)\right)}^{1/\alpha}f \circ \phi_t(s),\;\; t \in
\mathbb{Z}^d\,, \label{eqn4}
\end{eqnarray}
where, $f \in L^{\alpha}(S,\mu)$, $\{\phi_t\}_{t \in \mathbb{Z}^d}$
is a nonsingular $\mathbb{Z}^d$-action on $(S, \mu)$ and $\{c_t\}_{t
\in \mathbb{Z}^d}$ is a measurable cocycle for $\{\phi_t\}$ taking
values in $\{-1, +1\}$ i.e. each $c_t$ is a measurable map $c_t:S
\rightarrow \{-1, +1\}$ such that $\forall \, u,v \in \mathbb{Z}^d$
\[
c_{u+v}(s)=c_v(s) c_u\big(\phi_v(s)\big)\;\;\mbox{for $\mu$-a.a. }s
\in S.
\]
Conversely, if $\{f_t\}$ is of the form $(\ref{eqn4})$ then
$\{X_t\}$ defined by $(\ref{eqn3})$ is a stationary $S\alpha S$
random field. In particular, every {\it minimal} representation of the
process (see \cite{hardin:1982b}) turns out to be of the form
(\ref{eqn4}). \\

\noindent We will say that a stationary $S\alpha S$ random field
$\{X_t\}_{t\in \mathbb{Z}^d}$ is generated by a nonsingular
$\mathbb{Z}^d$-action $\{\phi_t\}$ on $(S, \mu)$ if it has a
integral representation of the form $(\ref{eqn4})$ satisfying
$(\ref{eqn5})$. With this terminology, we have the following extension
of Theorem 4.1 in \cite{rosinski:1995} to random fields.

\begin{propn} \label{propn2}
Suppose $\{X_t\}_{t\in \mathbb{Z}^d}$ is a stationary $S\alpha S$
random field generated by a nonsingular $\mathbb{Z}^d$-action
$\{\phi_t\}$ on $(S, \mu)$ and $\{f_t\}$ is given by $(\ref{eqn4})$.
Also let, $\mathcal{C}$ and $\mathcal{D}$ be the conservative and
dissipative parts of $\{\phi_t\}$. Then we have,
\begin{eqnarray*}
\mathcal{C}&=&\{s \in S:\sum_{t \in \mathbb{Z}^d}
|f_t(s)|^{\alpha}=\infty
\}\;\;\mbox{mod }\mu,\,and,\\
\mathcal{D}&=&\{s \in S:\sum_{t \in \mathbb{Z}^d}
|f_t(s)|^{\alpha}<\infty \}\;\;\mbox{mod }\mu\,.
\end{eqnarray*}
In particular, if a stationary $S\alpha S$ random field $\{X_t\}_{t \in
\mathbb{Z}^d}$ is generated by a conservative (dissipative, resp.)
$\mathbb{Z}^d$-action, then in any other integral representation of
$\{X_t\}$ of the form $(\ref{eqn4})$ satisfying $(\ref{eqn5})$, the
$\mathbb{Z}^d$-action must be conservative (dissipative, resp.).
Hence the classes of stationary $S\alpha S$ random fields generated
by conservative and dissipative actions are disjoint.
\end{propn}

\begin{proof}

Define $g$ as
\[
g(s)=\sum_{u \in {\mathbb{Z}}^d} \alpha_u \frac{d\mu \circ \phi_u}{d
\mu}(s)|f \circ \phi_u(s)|^{\alpha}
\]
$\mbox{where }\alpha_u > 0$ for all $u \in {\mathbb{Z}}^d$ and $\sum_{u
\in {\mathbb{Z}}^d} \alpha_u = 1.$ Clearly $g \in L^1$ and, by
$(\ref{eqn5})$, $g>0 \mbox{ a.e. }\mu$. Since
\[
\sum_{t \in \mathbb{Z}^d} \frac{d \mu \circ \phi_t}{d \mu}(s) g
\circ \phi_t(s) = \sum_{t \in \mathbb{Z}^d} \frac{d \mu \circ
\phi_t}{d \mu}(s) |f \circ \phi_t(s)|^{\alpha}=\sum_{t \in
\mathbb{Z}^d} |f_t(s)|^{\alpha}
\]
we can use Corollary $\ref{dualcor}$ to establish the first part of the
proposition, from which the second part of the proposition follows by
the same argument as in the one-dimensional case (see Theorem 4.1 in
\cite{rosinski:1995}).
\end{proof}

As in the one-dimensional case, it follows that the test described in
the previous proposition can be applied to any full support integral
representation of the process, not necessarily that of a specific
form.

\begin{cor} \label{cor_gen_cons_diss} The stationary $S\alpha S$
  random field $\{X_t\}_{t \in
\mathbb{Z}^d}$ is generated by a conservative (dissipative, resp.)
$\mathbb{Z}^d$-action if and only if for any (equivalently, some)
  integral representation
$(\ref{eqn3})$ of $\{X_t\}$ satisfying $(\ref{eqn5})$, the sum
\[
\sum_{t \in \mathbb{Z}^d} |f_t(s)|^{\alpha}
\]
is infinite (finite, resp) $\mu$-a.e. .
\end{cor}
 Proposition $\ref{propn2}$ also enables us to extend the
connection between the structure of stationary stable processes and
ergodic theory of nonsingular actions (given in \cite{rosinski:1995}) to
the case of stationary stable random fields. A decomposition of a
stable random field into three independent parts is available in
\cite{rosinski:2000}. A connection with the conservative-dissipative
decomposition is still missing in the case of random fields. Here we
provide the missing link. Recall that a stable random field
$\mathbf{X}$ is called a mixed moving average if it can be represented
in the form
\begin{eqnarray}
\mathbf{X} \eqdef \left\{\int_{W \times
{\mathbb{Z}}^d}f(v,t+s)\,M(dv,ds)\right\}_{t \in {\mathbb{Z}}^d}\,,
\label{mixed_moving_avg_defn}
\end{eqnarray}
where $f \in L^{\alpha}(W \times {\mathbb{Z}}^d, \nu \otimes l)$,
$l$ is the counting measure on ${\mathbb{Z}}^d$, $\nu$ is a
$\sigma$-finite measure on a standard Borel space $(W,
\mathcal{W})$, and the control measure $\mu$ of $M$ equals $\nu
\otimes l$ (see \cite{surgailis:rosinski:mandrekar:cambanis:1993} and
\cite{rosinski:2000}). The following result gives two equivalent
characterizations of stationary $S\alpha S$ random fields generated by
dissipative actions.

\begin{thm} \label{mixed_moving_avg_thm} Suppose $\{X_t\}_{t \in
{\mathbb{Z}}^d}$ is a stationary $S\alpha S$ random field. Then, the
  following are equivalent:
\begin{enumerate}
\item $\{X_t\}$ is generated by a dissipative $\mathbb{Z}^d$-action.

\item For any integral representation $\{f_t\}$ of $\{X_t\}$ we have,
\[
\sum_{t \in \mathbb{Z}^d} |f_t(s)|^{\alpha}<\infty \mbox{ for }
\mu\mbox{-a.a. }s.
\]
\item $\{X_t\}$ is a mixed moving average.

\end{enumerate}

\end{thm}

\begin{proof} $1$ and $2$ are equivalent by Corollary
$\ref{cor_gen_cons_diss}$, and, $2$ and $3$ are equivalent by Theorem $2.1$ of
\cite{rosinski:2000}.

\end{proof}

Theorem $\ref{mixed_moving_avg_thm}$ allows us to describe the
decomposition of a stationary $S \alpha S$ random field given in
Theorem 3.7 of \cite{rosinski:2000} in terms of the
ergodic-theoretical properties of  nonsingular
$\mathbb{Z}^d$-actions generating the field. The statement of the
following corollary is an extension of the one-dimensional
decomposition in Theorem 4.3 in \cite{rosinski:1995} to random
fields.

\begin{cor} \label{c:decomp}
A stationary $S \alpha S$ random field $\mathbf{X}$ has a unique in
law decomposition
\begin{eqnarray}
X_t&\eqdef& X^{\mathcal{C}}_t+X^{\mathcal{D}}_t \,,
\label{cons_diss_decomp}
\end{eqnarray}
where $\mathbf{X^{\mathcal{C}}}$ and $\mathbf{X^{\mathcal{D}}}$ are
two independent stationary $S\alpha S$ random fields such that
$\mathbf{X^{\mathcal{D}}}$ is a mixed moving average, and
$\mathbf{X^{\mathcal{C}}}$ is generated by a conservative action.
\end{cor}

As mentioned before, all of the structure results of this section
extend immediately to $G$-stationary random fields for countable
abelian groups $G$ more general than  $\mathbb{Z}^d$. The only place
where an additional argument is needed is the equivalence of parts 2
and 3 in Theorem \ref{mixed_moving_avg_thm}, with a $G$-mixed moving
average defined by
\[
\mathbf{X} \eqdef \left\{\int_{W \times
G}g(v,t+s)\,M(dv,ds)\right\}_{t \in G}\,,
\]
in notation parallel to \eqref{mixed_moving_avg_defn}. This
equivalence needs an extension of Theorem $2.1$ in
\cite{rosinski:2000} to general countable abelian groups. See
\cite{roy:2007} for details of this extension which does not require
any additional ideas to what is already in the original proof.

As in the one-dimensional case, it is possible to think of
stable random fields generated by conservative actions as having
longer memory than those  generated by dissipative actions, simply
because a conservative action ``keeps
coming back'', and so the same values of the random measure $M$
contribute to observations $X_t$ far separated in $t$. From this point
of view, the $\mathbb{Z}^d$-action $\{\phi_t\}$ is a parameter
(though highly infinite-dimensional) of the stationary $S \alpha S$
random field $\{X_t\}$ that determines, among others, the length of
its memory.


\section{Maxima of Stable Random Fields}\label{b_n_section}

\noindent The length of  memory of stable random fields is
manifested, in particular, in the rate of growth of its extreme
values. If $X_t$ is
generated by a conservative action, the extreme values tend to grow at a
slower rate because longer memory prevents erratic changes in $X_t$
even when $t$ becomes ``large''. This has been
formalized in \cite{samorodnitsky:2004a} for $d=1$, and it turns out to
be the case for  stable random fields as well.

For a stationary $S \alpha S$ random field $\{X_t\}_{t \in
\mathbb{Z}^d}$, we will study the partial maxima sequence
\begin{eqnarray}
M_n &:=& \max_{0 \leq t \leq (n-1)\mathbf{1}}|X_t|, \;\;\; n
=0,1,2,\ldots \label{maxm}
\end{eqnarray}
where $u = (u^{(1)}, u^{(2)}, \ldots, u^{(d)}) \leq v =(v^{(1)},
v^{(2)}, \ldots, v^{(d)})$ means $u^{(i)} \leq v^{(i)}$ for all
$i=1,2,\ldots,d$ and $\mathbf{1}=(1,1,\ldots,1)$. As in the
one-dimensional case, the asymptotic behavior of the maximum functional
$M_n$ is related to the deterministic sequence
\begin{eqnarray}
b_n&:=& \left(\int_S \max_{0 \leq t \leq (n-1)
\mathbf{1}}|f_t(s)|^{\alpha}\mu(ds)\right)^{1/\alpha}, \;\;\;n
=0,1,2,\ldots \,\,. \label{b_n}
\end{eqnarray}
Note that $b_n$ is completely determined by the process, and does
not depend on a particular integral representation (see Corollary
$4.4.6$ of \cite{samorodnitsky:taqqu:1994}). We are interested in the
features of this sequence that are related to the decomposition of a
stable random field in Corollary \ref{c:decomp}.  The next result
shows that the sequence $b_n$ grows at a slower rate for random
fields generated by a conservative action than for random fields
generated by a dissipative action.

\begin{propn} \label{b_n_rate_propn} Let $\{f_t\}$ be given by
$(\ref{eqn4})$. Assume that $(\ref{eqn5})$ holds.
 \begin{enumerate}
 \item If the action $\{\phi_t\}$ is conservative then:
   \begin{eqnarray}
   n^{-d/\alpha}b_n&\rightarrow&0\;\;\;\mbox{as }n\rightarrow \infty.
\label{b_n_rate_str_cons}
   \end{eqnarray}
 \item If the action $\{\phi_t\}$ is dissipative, and the random field is given in the mixed moving
 average form $(\ref{mixed_moving_avg_defn})$, then:
   \begin{eqnarray}
   \lim_{n \rightarrow \infty} n^{-d/\alpha}b_n&=&{\left(\int_W {(g(v))}^\alpha
\nu(dv) \right)}^{1/\alpha} \in (0,\infty),\label{b_n_rate_diss}
   \end{eqnarray}
 where
   \begin{eqnarray}
   g(v)= \sup_{s \in \mathbb{Z}^d} |f(v,s)| \;\;\; \mbox{for }v \in W.
\label{g_defn}
   \end{eqnarray}
 \end{enumerate}
\end{propn}

\begin{proof} 1. Firstly we observe that without loss of generality we
can assume
that $\mu$ is a probability measure. This is because if $\nu$ is a
probability measure equivalent to the $\sigma$-finite measure $\mu$
then instead of $(\ref{eqn3})$ we will use
\[
X_t \eqdef \int_{S} h_t(s)N(ds)
\]
where,
\[
h_t(s) = c_t(s){\left(\frac{d \nu \circ \phi_t}{d
\nu}(s)\right)}^{1/\alpha} h \circ \phi_t(s),\;\; t \in \mathbb{Z}^d
\]
where $h=f \big(\frac{d\mu}{d\nu} \big)^{1/\alpha} \in L^\alpha(S,
\nu)$ and $N$ is a $S \alpha S$ random measure on $S$ with
control measure $\nu$.\\

Since $\{b_n\}$ is an increasing sequence, it is enough to show
$(\ref{b_n_rate_str_cons})$ along the odd subsequence. By
stationarity of $\{X_t\}$, we need to check that
\[
a_n:=\frac{1}{(2n+1)^d} \int_S \max_{t \in J_n} |f_t(s)|^\alpha
\mu(ds) \rightarrow 0\,,
\]
where $J_n:=\{(i_1,i_2,\ldots,i_d):-n \leq i_1,i_2,\ldots,i_d \leq n
\}$. Let $g=|f|^\alpha $. Then $\|g\|:= \int_S g(s) \mu(ds) <
\infty$, and we have for $0<\epsilon<1$
\begin{eqnarray*}
a_n&=&\frac{1}{(2n+1)^d} \int_S \max_{t \in J_n} \hat{\phi}_t g(s)
\mu(ds)\\
&\leq&\frac{1}{(2n+1)^d}\bigg(\int_S \max_{t \in J_n}
\Big[\hat{\phi}_t g(s)I\Big(\hat{\phi}_t g(s)\leq \epsilon \sum_{u
\in
J_n} \hat{\phi}_u g(s)\Big)\Big] \mu(ds)\\
&&\;\;\;\;\;\;\;\;\;\;\;\;\;\;\;\;\;+\int_S \max_{t \in J_n}
\Big[\hat{\phi}_t g(s)I\Big(\hat{\phi}_t g(s)> \epsilon \sum_{u \in
J_n} \hat{\phi}_u g(s)\Big)\Big] \mu(ds)\bigg)\\
&&=a^{(1)}_n+a^{(2)}_n\,.
\end{eqnarray*}
Clearly,
\begin{equation}
a^{(1)}_n \leq \frac{\epsilon}{(2n+1)^d} \sum_{u \in J_n} \int_S
\hat{\phi}_u g(s) \mu(ds) =\epsilon \|g\| \;,
\label{bound_on_a^{(1)}_n}
\end{equation}
and
\begin{equation}
a^{(2)}_n \leq  \frac{1}{(2n+1)^d} \sum_{t \in J_n} \int_S
\hat{\phi}_t g(s)I_{A_{t,n}}(s) \mu(ds)\,,
\label{bound_on_a^{(2)}_n}
\end{equation}
where $A_{t,n}=\{s:\, \hat{\phi}_t g(s) > \epsilon \sum_{u \in J_n}
\hat{\phi}_u g(s)\}\,,\,n \geq 1,\, t \in J_n\,$. Notice that for all
$n \geq 1$, and, for all $t \in
J_n$\,,
\begin{equation}
\int_S \hat{\phi}_t g(s)I_{A_{t,n}}(s) \mu(ds)=\int_S g(s)
I_{\phi_t^{-1}(A_{t,n})}(s) \mu(ds)\,.
\label{simplification_of_integral}
\end{equation}
The following is the most important step of this proof: if we define
\[
U_n:=\{(t_1,t_2,\ldots,t_d): -n+[\sqrt{n}] \leq t_1,t_2,\ldots,t_d
\leq n-[\sqrt{n}]\}
\]
then we have,
\begin{equation}
\lim_{n \rightarrow \infty} \max_{t \in
U_n}\;\mu(\phi_t^{-1}(A_{t,n}))=0\;. \label{limit_of_max_mu}
\end{equation}
To prove $(\ref{limit_of_max_mu})$ observe that for all $t \in U_n$
\begin{eqnarray*}
&&\phi_t^{-1}(A_{t,n})\\
&&\;\;\;\;=\Big\{\phi_{-t}(s):\,g\circ\phi_{-t}(s)\frac{d\mu\circ\phi_{-t}}{d\mu}(s)>\epsilon
\sum_{u \in J_n} g\circ\phi_{-u}(s)\frac{d\mu\circ\phi_{-u}}{d\mu}(s)\Big\}\\
&&\;\;\;\;=\Big\{s:\,g(s)>\epsilon \sum_{u \in J_n} g\circ\phi_{u+t}(s)\frac{d\mu\circ\phi_{u+t}}{d\mu}(s)\Big\}\\
&&\;\;\;\;\subseteq \Big\{s:\,g(s)>\epsilon \sum_{\tau \in
J_{[\sqrt{n}]}}
g\circ\phi_\tau(s)\frac{d\mu\circ\phi_\tau}{d\mu}(s)\Big\}\,.
\end{eqnarray*}
The last inclusion holds because $J_{[\sqrt{n}]} \subseteq t+J_n$. Hence, for any $M>0$
\begin{eqnarray*}
\max_{t \in U_n}\;\mu(\phi_t^{-1}(A_{t,n}))&\leq& \mu\{s:\,g(s)>\epsilon M\}+\mu\bigg(\sum_{t \in J_{[\sqrt{n}]}} g\circ\phi_t\frac{d\mu\circ\phi_t}{d\mu}\leq M\bigg)\\
&\leq& \frac{\|g\|}{\epsilon M}+\mu\Big(\sum_{t \in J_{[\sqrt{n}]}}
|f_t|^{\alpha}\leq M\Big)\,.
\end{eqnarray*}
Now $(\ref{limit_of_max_mu})$ follows by first using Proposition
$\ref{propn2}$ with a fixed $M$ and
 then letting $M\to\infty$. \\

\noindent From $(\ref{simplification_of_integral})$ and
$(\ref{limit_of_max_mu})$ it follows that
\begin{eqnarray}
&&\frac{1}{(2n+1)^d} \sum_{t \in U_n} \int_S \hat{\phi}_t
  g(s)I_{A_{t,n}}(s) \mu(ds) \nonumber \\
&&\;\;\;\;\;\;\;\;\;\;\;\;\;\;\;\;\;\;\;\;\;\;\;\;\;\;\;\;=\frac{1}{(2n+1)^d}
\sum_{t \in U_n}\int_{\phi_t^{-1}(A_{t,n})} g(s) \mu(ds) \rightarrow
0\,. \label{bound_on_part_of_a^{(2)}_n}
\end{eqnarray}

\noindent If we define $V_n=J_n-U_n$, then
$$
\frac{1}{(2n+1)^d} \sum_{t \in V_n} \int_S \hat{\phi}_t
g(s)I_{A_{t,n}}(s) \mu(ds)
\leq \frac{1}{(2n+1)^d} \sum_{t \in V_n} \int_S
\hat{\phi}_t g(s)\mu(ds) \to 0\,.
$$
Then using $(\ref{bound_on_a^{(2)}_n})$ and
$(\ref{bound_on_part_of_a^{(2)}_n})$ we see that $a^{(2)}_n\to 0$ as
$n\to\infty$. Therefore we get,
\[
\limsup a_n \leq \limsup a^{(1)}_n + \limsup a^{(2)}_n \leq \epsilon
\|g\|\,,
\]
and, since $\epsilon >0$ is arbitrary, the result follows.\\

\noindent 2. The argument here is similar to that used in the
one-dimensional case in Theorem $3.1$ of \cite{samorodnitsky:2004a}.
One uses a direct computation to check the claim in the case where
$f$ has compact support, that is
\[
f(v,s)I_{W \times [-m\mathbf{1},m\mathbf{1}]^{c}}(v,s) \equiv 0
\mbox{ for some }m=1,2,\ldots
\]
where $[u,v]:=\{t \in \mathbb{Z}^d: u \leq t \leq v\}$. The
proof in the general case follows then by approximating a general kernel
$f$ by a kernel with a compact support.

\end{proof}

\begin{remark} \label{rk:small.b.n.F}
\textnormal{The statement of the first part of the proposition clearly
  extends  to $G$-stationary random fields for any free abelian
  group $G$ of rank $d$, since the same is true for Proposition
$\ref{propn2}$. See the discussion after Corollary \ref{c:decomp}.}
\end{remark}

We are now ready to investigate the rate of growth of the sequence
$\{M_n\}$ of partial maxima of a stationary symmetric
$\alpha$-stable random field, $0 < \alpha < 2$. We will see that if
such a random field has a nonzero component $X^{\mathcal{D}}$ in
$(\ref{cons_diss_decomp})$ generated by a dissipative action, then
the partial maxima grow at the rate $n^{d/\alpha}$, while if the
random field is generated by a conservative action, then the partial
maxima grow at a slower rate. As we will see in the sequel, the
actual rate of growth of the sequence $\{M_n\}$ in the conservative
case, depends on a number of factors. The dependence on the group
theoretical properties of the action is very prominent. We start
with the following result, which extends  Theorem $4.1$ of
\cite{samorodnitsky:2004a} to $d>1$. It is based on Proposition
\ref{b_n_rate_propn}, and the argument is parallel to the
one-dimensional case.

\begin{thm} \label{M_n_rate_thm}
Let $\mathbf{X}=\{X_t\}_{t \in \mathbb{Z}^d}$ be a stationary
$S\alpha S$ random field, with $0 < \alpha < 2$, integral
representation $(\ref{eqn3})$, and
functions $\{f_t\}$ given by $(\ref{eqn4})$. \\

\noindent 1. Suppose that $\mathbf{X}$ is not generated by a
conservative action (i.e. the component $X^{\mathcal{D}}$ in
$(\ref{cons_diss_decomp})$ generated by the dissipative part is
nonzero). Then
\begin{eqnarray}
\frac{1}{n^{d/\alpha}}M_n  \Rightarrow  C^{\,1/\alpha}_\alpha K_X
Z_\alpha \label{M_n_rate_noncons}
\end{eqnarray}
as $n \rightarrow \infty$, where
\[
K_X ={\left(\int_W {(g(v))}^\alpha \nu(dv) \right)}^{1/\alpha}
\]
and $g$ is given by $(\ref{g_defn})$ for any representation of
$X^{\mathcal{D}}$ in the mixed moving average form
$(\ref{mixed_moving_avg_defn})$, $C_\alpha$ is the stable tail
constant (see (1.2.9) in \cite{samorodnitsky:taqqu:1994})
and $Z_\alpha$ is the standard Frechet-type extreme value random
variable with the distribution
\[
P(Z_\alpha \leq z)=e^{-z^{-\alpha}},\;\;\;\;z>0.
\]

\noindent 2. Suppose that $\mathbf{X}$ is generated by a
conservative $\mathbb{Z}^d$-action. Then
\begin{eqnarray}
\frac{1}{n^{d/\alpha}}M_n \probconv 0 \label{M_n_rate_str_cons}
\end{eqnarray}
as $n \rightarrow \infty$. Furthermore, with $b_n$ given by
$(\ref{b_n})$,
\begin{eqnarray}
&&\left\{\frac{1}{c_n}\,M_n \right\} \mbox{ is not tight for any
positive sequence }c_n=o(b_n), \label{nontightness}
\end{eqnarray}
while
\begin{eqnarray}
&&\left\{\frac{1}{b_n \zeta_n}\, M_n \right \} \mbox{ is tight,
where }\zeta_n = \left\{
  \begin{array}{ll}
  1,                            & \mbox{ if }\; 0 < \alpha < 1,\\
  L_n,                          & \mbox{ if }\;\;\; \alpha = 1,\\
  (\log{n})^{1/\alpha^{\prime}},& \mbox{ if }\; 1 < \alpha < 2\,,
  \end{array}
 \right. \label{tightness}
\end{eqnarray}
where $L_n:=\max(1, \log{\log{n}})$, and for $\alpha>1$,
$\alpha^{\prime}$ is such that $1/\alpha + 1/\alpha^{\prime}=1$. \\

\noindent If, for some $\theta >0$ and $c >0$,
\begin{eqnarray}
b_n  \geq  cn^\theta \;\;\;\;\text{ for all $n \geq 1,$} \label{inequality}
\end{eqnarray}
then $(\ref{tightness})$ holds with $\zeta_n \equiv 1 $ for all
$0 < \alpha < 2$.\\

\noindent Finally, for $n=1,2,\ldots$, let $\eta_n$ be a
probability measure on $(S, \mathcal{S})$ with
\begin{eqnarray}
\frac{d \eta_n}{d \mu}(s)=b_n^{-\alpha} \max_{0 \leq t \leq
(n-1)\mathbf{1}} |f_t(s)|^\alpha,&& s \in S,
\label{density_of_eta_n}
\end{eqnarray}
and let $U_j^{(n)},\, j=1,2$ be independent $S$-valued random
variables with common law $\eta_n$. Suppose that
$(\ref{inequality})$ holds and for any $\epsilon
> 0$,
\begin{eqnarray}
&&P\bigg(\mbox{for some }t \in [0, (n-1)\mathbf{1}],\nonumber\\
&&\;\;\;\;\; \frac{|f_t(U_j^{(n)})|}{\max_{0 \leq u \leq (n-1)
\mathbf{1}} |f_u(U_j^{(n)})|} > \epsilon,\; j=1,2 \bigg)\rightarrow
0 \label{suff_condn}
\end{eqnarray}
as $n \rightarrow \infty$. Then
\begin{eqnarray}
\frac{1}{b_n}M_n  \Rightarrow  C^{\,1/\alpha}_\alpha Z_\alpha
\label{M_n_rate_under_suff_condn}
\end{eqnarray}
as $n \rightarrow \infty$.
\end{thm}

\begin{remark} \label{suff_condn_remark}
\textnormal{An easily verifiable sufficient condition for
$(\ref{suff_condn})$ is
\begin{eqnarray}
\lim_{n \rightarrow \infty} \frac{b_n}{n^{d/2 \alpha}} = \infty.
\label{suff_condn_2}
\end{eqnarray}
Alternatively, $(\ref{suff_condn})$ holds if we assume that $\mu$ is
a finite measure, $\{\phi_t\}$ is measure preserving, the sequence
$\{b_n^{-\alpha}\max_{0 \leq t \leq
(n-1)\mathbf{1}}|f_t(s)|^\alpha\}$, $t\in \mathbb{Z}^d$ is uniformly
integrable with respect to $\mu$ and, for every $\epsilon >0$
\begin{equation}
\lim_{n \rightarrow \infty} n^{d/2}\mu\{s \in S:\,|f(s)|>\epsilon
b_n\}=0\,. \label{suff_condn_3}
\end{equation}
The arguments are the same as in the case $d=1$.}
\end{remark}


\section{Connections with Group Theory} \label{sec_group_theory}

When the underlying action is not conservative Theorem
$\ref{M_n_rate_thm}$ yields the exact rate of growth of $M_n$. For
conservative actions, however, the actual rate of growth of the
partial maxima depends on further properties of the action. In this
section we investigate the effect of the group theoretic structure
of the action on the  rate of growth of the partial maximum. We
start with introducing the
appropriate notation.\\

Consider $A:=\{\phi_t:\,t \in \mathbb{Z}^d\}$ as a subgroup of the
group of invertible nonsingular transformations on $(S,\mu)$ and
define a group homomorphism
\[
\Phi:\mathbb{Z}^d \rightarrow A
\]
by $\Phi(t)=\phi_t$ for all $t \in \mathbb{Z}^d$. Let
$K:=Ker(\Phi)=\{t \in \mathbb{Z}^d:\,\phi_t = 1_S\}$, where $1_S$
denote the identity map on $S$. Then $K$ is a free abelian group and
by first isomorphism theorem of groups (see, for example,
\cite{lang:2002}) we have,
\[
A \simeq \mathbb{Z}^d/K \,.
\]
Hence by Theorem $8.5$ in Chapter I of \cite{lang:2002} we get,
\[
A=\bar{F} \oplus \bar{N}
\]
where $\bar{F}$ is a free abelian group and $\bar{N}$ is a finite
group. Assume $rank(\bar{F})=p \geq 1$ and $|\bar{N}|=l$. Since
$\bar{F}$ is free, there exists an injective group homomorphism
\[
\Psi: \bar{F} \rightarrow \mathbb{Z}^d
\]
such that $\Phi \circ \Psi = 1_{\bar{F}}$. Let $F=\Psi(\bar{F})$.
Then $F$ is a free subgroup of $\mathbb{Z}^d$ of rank $p$.\\

The rank $p$ is the effective dimension of the random field, giving
more precise information on the rate of growth of the partial
maximum than the nominal dimension $d$. We start with showing that
this is true for the sequence $\{b_n\}$ in \eqref{b_n}.

\begin{propn}\label{propn_b_n_Grate} Let $\{f_t\}$ be given by
$(\ref{eqn4})$. Assume that $(\ref{eqn5})$ holds. Then we have the
  following:

\begin{enumerate}

 \item If $\{\phi_t\}_{t \in F}$ is conservative then
 \begin{equation}
 n^{-p/\alpha}b_n \rightarrow 0 \,.\label{b_n_rate_Gcons}
 \end{equation}

 \item If $\{\phi_t\}_{t \in F}$ is dissipative then
 \begin{equation}
 n^{-p/\alpha}b_n \rightarrow a \label{b_n_rate_Gdiss}
 \end{equation}
 for some $a \in (0, \infty)$.

\end{enumerate}
\end{propn}

\begin{proof} 1. It is easy to check that $F \cap K = \{0\}$ and hence the sum $F+K$ is
direct. Suppose $G=F \oplus K$. Using group isomorphism theorems we
have,
\[
\mathbb{Z}^d/G \simeq  (\mathbb{Z}^d/K)/(F \oplus K/K) \simeq
A/\bar{F} \simeq  \bar{N}\,.
\]
Assume that $x_1+G,\, x_2+G,\,\ldots\,,x_l+G$ are all the cosets of
$G$ in $\mathbb{Z}^d$. Let $rank(K)=q$. Choose a basis
$\{u_1,u_2,\ldots,u_p\}$ of $F$ and a basis $\{v_1,v_2,\ldots,v_q\}$
of $K$. We need the following\\

\begin{lemma} \label{lemma_inclusions} There are positive integers
  $c$, $d$, and, $N$ such that for every $n\geq N$
\begin{equation}
\bigcup_{k=1}^{l} (x_k+G_{[n/d]}) \subseteq [-n\mathbf{1},
n\mathbf{1}] \subseteq \bigcup_{k=1}^{l} (x_k+G_{cn})
\label{inclusions}
\end{equation}
where for $m\geq 1$
\[
G_m:=\{\sum_{i=1}^{p} \alpha_i u_i + \sum_{j=1}^q \beta_j v_j:\;
|\alpha_i|, |\beta_j| \leq m \; \text{\rm for all}\; i,j\}\,.
\]
\end{lemma}

\begin{proof} Let $r=p+q$. For ease of notation we define
\[
w_i = \left \{
            \begin{array}{lc}
            u_i \;\;\;\;& 1 \leq i \leq p\,, \\
            v_{i-p} \;\;\; & \; p+1 \leq i \leq r\,. \\
            \end{array}
      \right .
\]
Then $\{w_1,w_2,\ldots,w_r\}$ is a basis for $G$. The first
inclusion in $(\ref{inclusions})$ is obvious. To establish the
second inclusion we first prove\\

\noindent \textbf{\textit{Step 1.}} There is an integer $c^{\prime}
\geq 1$ such that
\[
[-n\mathbf{1}, n\mathbf{1}] \cap G \subseteq G_{c^{\prime}n}
\;\;\;\mbox{ for all }n \geq 1\,.
\]

\begin{proof}[Proof of Step 1]
\renewcommand{\qedsymbol}{}
Take $y \in [-n\mathbf{1}, n\mathbf{1}] \cap G$. Then, $y=\eta_1 w_1
+ \eta_2 w_2 + \cdots + \eta_r w_r$ for some $\eta_1, \eta_2,
\ldots, \eta_r \in \mathbb{Z}$\,. We have to show $| \eta_i | \leq
c^{\prime} n$ for all $1 \leq i \leq r$ for some $c^{\prime} \geq 1$
that does not depend on $n$. Let $\widetilde{\eta}^{T}:=(\eta_1,
\eta_2, \ldots,\eta_r) \in \mathbb{Z}^r$. Then,
\begin{equation}
y=W\widetilde{\eta} \label{form_of_y}
\end{equation}
where, $W$ is the $d \times r$ matrix with $w_i$ as the $i^{th}$
column. The columns of $W$ are linearly independent over
$\mathbb{Z}$ and hence over $\mathbb{R}$. Hence there is a $r \times
d$ matrix $Z$ such that
\[
ZW=I
\]
where $I$ is the identity matrix of order $r$. Hence from
$(\ref{form_of_y})$ we have,
\[
\widetilde{\eta}=Zy\,.
\]
For all $1 \leq i \leq r$ we get,
\begin{equation*}
|\eta_i|\, \leq \,\|\widetilde{\eta}\|\, \leq \, \|Z\|\|y\|\, \leq
\, \|Z\|n\sqrt{d}\, \leq \, c^{\prime} n
\end{equation*}
where, $c^{\prime}=\bigl[\|Z\|\sqrt{d}\bigr]+1$.
This proves Step $1$. \\
\end{proof}
\noindent \textbf{\textit{Step 2.}} Let
\[
M=\max_{1 \leq k \leq l} \|x_k\|_{\infty}+1
\]
where $\|\cdot\|_{\infty}$ denotes the sup-norm on $R^d$, and
$c=c^{\prime} M$. Then for all $n \geq 1$ we have,
\[
[-n\mathbf{1}, n\mathbf{1}] \subseteq \bigcup_{k=1}^{l}
(x_k+G_{cn})\;.
\]

\begin{proof}[Proof of Step 2]
\renewcommand{\qedsymbol}{}
Take $y \in [-n\mathbf{1}, n\mathbf{1}]$. Then $y \in x_{k_0}+G$ for
some $1 \leq k_0 \leq l$. Clearly, $y^{\prime}:=y-x_{k_o} \in
[-(n+M-1)\mathbf{1}, (n+M-1)\mathbf{1}] \cap G$. By Step $1$,
$y^{\prime} \in G_{c^{\prime}(n+M-1)} \subseteq G_{cn}$, and hence,
$y \in x_{k_0}+G_{cn} \subseteq \bigcup_{k=1}^{l} (x_k+G_{cn})$,
proving Step 2 and the lemma.
\end{proof}
\end{proof}

For $k=1,\ldots ,l$ let
\[
g_k=f\circ\phi_{x_k}
\bigg(\frac{d\mu\circ\phi_{x_k}}{d\mu}\bigg)^{1/\alpha}\,.
\]
Then for $t=x_k+\sum_{i=1}^{p}
\alpha_i u_i + \sum_{j=1}^q \beta_j v_j$ we have
\begin{equation}
|f_t(s)|=|g_k \circ \phi_{\sum_{i=1}^{p} \alpha_i
u_i}(s)|\bigg(\frac{d\mu\circ \phi_{\sum_{i=1}^{p} \alpha_i
u_i}}{d\mu}(s)\bigg)^{1/\alpha} \label{form_of_modf_t} \,.
\end{equation}

By stationarity,  Lemma
$\ref{lemma_inclusions}$ and $(\ref{form_of_modf_t})$ we have, for
all $n \geq N$,
\begin{eqnarray*}
b_n^\alpha&\leq & b_{2n+1}^{\alpha}=\int_S \max_{-n\mathbf{1} \leq t
  \leq n\mathbf{1}}|f_t(s)|^{\alpha}\mu(ds)\\
&\leq& \int_S \max_{1 \leq k \leq l}\, \max_{|\alpha_i| \leq cn}
\bigg(|g_k \circ \phi_{\sum_{i=1}^{p} \alpha_i u_i}(s)|^\alpha
\frac{d\mu\circ \phi_{\sum_{i=1}^{p} \alpha_i u_i}}{d\mu}(s)\bigg)
\mu(ds) \\
&\leq & \sum_{k=1}^l \int_S \max_{|\alpha_i| \leq cn} \bigg(|g_k
\circ \phi_{\sum_{i=1}^{p} \alpha_i u_i}(s)|^\alpha \frac{d\mu\circ
\phi_{\sum_{i=1}^{p} \alpha_i u_i}}{d\mu}(s)\bigg) \mu(ds) \\
&=& o(n^{p})\,.
\end{eqnarray*}
The last step follows from Proposition $\ref{b_n_rate_propn}$ and
Remark \ref{rk:small.b.n.F}. \\

\noindent 2. Proof of this part is similar to the proof of Theorem
$3.1$ in \cite{samorodnitsky:2004a}. We start this proof with the
following combinatorial fact:

\begin{lemma} \label{lemma_density_of_F} For $n\geq 1$
and $k=1,2,\ldots,l$, let
\[
F_{k,n}=\bigl\{u \in x_k+F: \mbox{ there exists } v \in K \mbox{
such that } u+v \in [-n\mathbf{1}, n\mathbf{1}]\bigr\}\,.
\]
Then there is a positive real number $\mathcal{V}$ such that for all
$k=1,2,\ldots,l$,
\begin{equation}
\lim_{n \to \infty} \frac{|F_{k,n}|}{n^p}=\mathcal{V}\,.
\label{limit_density}
\end{equation}
Here $|A|$ stands for the cardinality of a set $A$.
\end{lemma}

\begin{proof} One of $F_{k,n}$ is the set
\[
F_n=\bigl\{y \in F: \mbox{ there exists } v \in K \mbox{ such that }
y+v \in [-n\mathbf{1}, n\mathbf{1}]\bigr\}\,.
\]
Firstly, we will show
\begin{equation}
\lim_{n \to \infty} \frac{|F_n|}{n^p}=\mathcal{V}
\label{limit_density_F}
\end{equation}
for some $\mathcal{V} > 0$. To show this let $W$ be the matrix used
in the proof of Lemma $\ref{lemma_inclusions}$. We can partition $W$
into two submatrices as follows:
\[
W=[\,U\,|\,V\,]
\]
where, $U$ is the $d \times p$ matrix whose $i^{th}$ column is $u_i$
and, $V$ is the $d \times q$ matrix whose $j^{th}$ column is $v_j$.
Since the columns of $U$ are linearly independent over $\mathbb{Z}$,
we have,
\[
|F_n|=|\{\alpha \in \mathbb{Z}^p:\mbox{there exists }\beta \in
\mathbb{Z}^q \mbox{ such that }\|U\alpha+V\beta\|_\infty \leq
n\}|\,.
\]
Let $P:=\{x \in \mathbb{R}^r:\|Wx\|_\infty \leq 1\}$ and
$\pi:\mathbb{R}^r \to \mathbb{R}^p$ denote the projection map on the
first $p$ coordinates:
\[
\pi(x_1,x_2,\ldots,x_r)=(x_1,x_2,\ldots,x_p)\,.
\]
Then we have,
\[
\frac{|F_n|}{n^p}=\frac{|\pi(\mathbb{Z}^r \cap nP)|}{n^p}=:a_n\,.
\]
Let
\[
b_n:=\frac{|\mathbb{Z}^p \cap n\pi(P)|}{n^p}\,.
\]
Clearly, $a_n \leq b_n$. Since $P$ is a rational polytope (i.e. a
polytope whose vertices have rational coordinates) so is $\pi(P)$.
 Hence, by Theorem $1$ of \cite{deloera:2005}, it follows that
\begin{equation}
\limsup_{n\to \infty} a_n \leq \lim_{n \to \infty} b_n = \mathcal{V}
\label{inequality_limsup}
\end{equation}
where $\mathcal{V}=Volume(\pi(P))$, the $p$-dimensional volume of
 $\pi(P)$. This volume is positive since the latter set, obviously,
 contains a small ball centered at the origin. For the other
 inequality we let
\[
P_m:=\left \{x \in \mathbb{R}^r:\|Wx\|_\infty \leq
1-\frac{\|W\|_\infty}{m}\right \}
\]
where $\|W\|_\infty:=\sup_{x \neq 0}
\frac{\|Wx\|_\infty}{\|x\|_\infty} \in \mathbb{Z}$ since $W$ is a
matrix with integer entries. Hence for all $m>\|W\|_\infty$, $P_m$
is a rational polytope of dimension $r$. Also, $P_m \uparrow P$. Now
fix $m > \|W\|_\infty$. Observe that
\begin{equation}
\left\{y \in \mathbb{R}^r: \|y-x\|_\infty \leq \frac{1}{m}\right\}
\subseteq P\;\;\;\mbox{ for all }x \in P_m\,. \nonumber
\end{equation}
Hence, it follows that for all $n>m$,
\[
\pi\left(\frac{1}{n}\mathbb{Z}^r \cap P\right) \supseteq
\frac{1}{n}\mathbb{Z}^p \cap \pi(P_m)\,,
\]
which, along with Theorem $1$ of \cite{deloera:2005}, implies
\begin{equation}
\liminf_{n\to \infty} a_n \geq \lim_{n\to \infty}\frac{|\mathbb{Z}^p
\cap n\pi(P_m)|}{n^p}=\mathcal{V}_m \label{inequality_liminf}
\end{equation}
where $\mathcal{V}_m=Volume(\pi(P_m))$ is, once again, the
$p$-dimensional volume. Since $P_m \uparrow P$,
it follows that $\mathcal{V}_m \uparrow \mathcal{V}$. Hence $(\ref{limit_density_F})$
follows from $(\ref{inequality_limsup})$ and $(\ref{inequality_liminf})$.\\

Now fix $k=1,2,\ldots,l$ and let $M=\|x_k\|$. Observe that for all
$n > M$,
\[
|F_{n-M}| \leq |F_{k,n}| \leq |F_{n+M}|\,.
\]
Hence $(\ref{limit_density})$ follows from
$(\ref{limit_density_F})$.

\end{proof}

We now return to the proof of the second part of the proposition.
We give a group structure to
\begin{equation}
H:=\bigcup_{k=1}^l (x_k + F) \label{defn_of_H}
\end{equation}
as follows.
For all $u_1, u_2 \in H$, there exists unique $u \in H$ such that
$(u_1+u_2)-u \in K$. We define this $u$ to be $u_1 \oplus u_2$. It
is not hard to check that $(H,\oplus)$ is a countable abelian group.
In fact, $H \simeq \mathbb{Z}^d/K$. We can define a nonsingular
group action $\{\psi_u\}$ of $H$ on $S$ as
\[
\psi_u = \phi_u \;\;\;\mbox{ for all }u \in H\,.
\]
Notice that if $h \in L^1(S, \mu), h>0$, then, since
$(\ref{defn_of_H})$ is a disjoint union,
\begin{equation}
\sum_{u \in H} \frac{d\mu \circ \psi_u}{d\mu} h \circ \psi_u =
\sum_{t \in F} \frac{d\mu \circ \phi_t}{d\mu} \tilde{h} \circ
\phi_t\,, \label{eqn_diss_of_psi}
\end{equation}
where,
\[
\tilde{h}=\sum_{k=1}^l \frac{d\mu \circ \phi_{x_k}}{d\mu}h \circ
\phi_{x_k}\,.
\]
Clearly $\tilde{h} \in L^1(S, \mu)$ and $\tilde{h}
> 0$. Hence using Corollary $\ref{dualcor}$ and the dissipativity of
$\{\phi_t\}_{t\in F}$, we see that the second sum in
\eqref{eqn_diss_of_psi} is finite almost everywhere. Another appeal
to Corollary $\ref{dualcor}$ shows that $\{\psi_u\}_{u \in H}$ is a
dissipative group action.

Define a random field $\{Y_u\}_{u\in H}$ as
\begin{equation} \label{e:H.field}
Y_u=\int_S \tilde{f}_u(s) M(ds),\;\;\;u\in H,
\end{equation}
where,
\[
\tilde{f}_u= f \circ \psi_u \left(\frac{d\mu \circ
\psi_u}{d\mu}\right)^{1/\alpha}\;\;\;u \in H.
\]
Clearly $\{Y_u\}_{u\in H}$ is an $H$-stationary $S\alpha S$ random field
generated by the dissipative action $\{\psi_u\}_{u\in H}$. Hence there
is a standard Borel space $(W,\mathcal{W})$ with a $\sigma$-finite
measure $\nu$ on it such that
$$
Y_u \eqdef \int_{W \times H} g(w,u \oplus s)\,N(dw,ds)\;\;\; u \in
H,
$$
for some $g \in L^{\alpha}(W \times H, \nu \otimes \tau)$, where
$\tau$ is the counting measure on $H$, and, $N$ is a $S\alpha S$
random measure on $W \times H$ with control measure $\nu \otimes
\tau$ (see the discussion following Corollary \ref{c:decomp}.)

Let, for all $w \in W$,
\begin{equation}
g^{\ast}(w):=\sup_{u \in H}|g(w,u)|\,.\label{defn_of_g^star}
\end{equation}
Then, clearly, $g^{\ast} \in L^\alpha(W,\nu)$. We will show that
(\ref{b_n_rate_Gdiss}) holds with
\begin{equation}
a:=\left(\frac{\mathcal{V}l}{2^p}\int_W (g^{\ast}(w))^\alpha
d\nu(w)\right)^{1/\alpha} \in (0,\infty)\,. \label{form_of_a}
\end{equation}
Since $b_n$ is an increasing sequence, it is enough to show
\begin{equation}
\lim_{n\to \infty}
\frac{b_{2n+1}}{{(2n+1)}^{p/\alpha}}=a\,.\label{b_n_rate_Gdiss_odd_subseq}
\end{equation}
Let $H_n:=\bigcup_{k=1}^l F_{k,n}$. Then by stationarity of
$\{X_t\}_{t \in \mathbb{Z}^d}$ we have, for all $n \geq 1$,
\begin{eqnarray}
b_{2n+1}^\alpha &=& \int_S \max_{-n\mathbf{1} \leq t \leq n\mathbf{1}}|f_t(s)|^\alpha \mu(ds) \nonumber \\
                &=& \int_S \max_{u \in H_n}|\tilde{f}_u(s)|^\alpha \mu(ds) \nonumber \\
                &=& \sum_{s \in H}\int_W \max_{u \in H_n} |g(w, s \oplus u)|^\alpha \nu(dw)\,. \label{form_of_b_n}
\end{eqnarray}
The last equality follows from Corollary $4.4.6$ of \cite{samorodnitsky:taqqu:1994}. We define a map $N:H \to \{0,1,\ldots\}$ as,
\[
N(u):=\min\{\|u+v\|_ \infty: v \in K\}\,.
\]
Clearly, for all $u \in H$,
\begin{equation}
N(u^{-1})=N(u)\,, \label{symmetry_of_N}
\end{equation}
where, $u^{-1}$ is the inverse of $u$ in $H$. Also, $N(\cdot)$
 satisfies the following ``triangle inequality": for all $u_1, u_2
\in H$,
\begin{equation}
N(u_1 \oplus u_2) \leq N(u_1) + N(u_2)\,.
\label{triangle_inequality_of_N}
\end{equation}
Observe that $H_n=\{u \in H: N(u) \leq n\}$. From Lemma
$\ref{lemma_inclusions}$ we have $H_n$'s are finite and Lemma
$\ref{lemma_density_of_F}$ yields
\begin{equation}
|H_n| \sim \mathcal{V}ln^p\,.\label{size_of_H_n}
\end{equation}
Also, clearly, $H_n \uparrow H$. As in the proof of Theorem $3.1$ of
\cite{samorodnitsky:2004a}, we first assume $g$ has compact support,
i.e. $g(w,u)I_{W \times H_m^{c}}(w,u)=0$ for some $m \geq 1$. Then
using $(\ref{symmetry_of_N})$ and
$(\ref{triangle_inequality_of_N})$, the expression in $(\ref{form_of_b_n})$
becomes
\begin{eqnarray*}
b_{2n+1}^\alpha&=& \sum_{s \in H_{n+m}}\int_W \max_{u \in H_n} |g(w, s \oplus u)|^\alpha \nu(dw)\\
                 &=& \sum_{s \in H_{n-m}}\int_W \max_{u \in H_n} |g(w, s \oplus u)|^\alpha \nu(dw) \\
                 &&\;\;\;\;\;\;\;+ \sum_{s \in H_{n+m}\cap H_{n-m}^{c}}\int_W \max_{u \in H_n} |g(w, s \oplus u)|^\alpha \nu(dw)=:A_n + B_n
\end{eqnarray*}
for all $n > m$. Using $(\ref{symmetry_of_N})$ and
$(\ref{triangle_inequality_of_N})$ once again, we have, for all $s \in H_{n-m}$,
\[
\max_{u\in H_n} |g(w,s\oplus u)|=g^{\ast}(w)\,.
\]
Hence, using $(\ref{size_of_H_n})$, we get,
\[
A_n=|H_{n-m}|\int_W (g^{\ast}(w))^\alpha \nu(dw) \sim a^\alpha (2n+1)^p\,,
\]
while
\[
B_n \leq \big(|H_{n+m}|-|H_{n-m}|\big)\int_W (g^{\ast}(w))^\alpha \nu(dw)=o(n^p)\,.
\]
Hence, $(\ref{b_n_rate_Gdiss_odd_subseq})$ follows for $g$ having compact support. The
proof in the general case follows by approximating a general kernel
$g$ by a kernel with a compact support as done in the proof of Theorem
$3.1$ in \cite{samorodnitsky:2004a}. This completes the proof of the
proposition.
\end{proof}

The following result sharpens the the description of the
asymptotic behaviour of the partial maxima of a random field given in
Theorem \ref{M_n_rate_thm}. It reduces to the latter result if
$K=0$.

\begin{thm}\label{thm_M_n_Grate} Let $\mathbf{X}=\{X_t\}_{t \in
    \mathbb{Z}^d}$ be a stationary
$S\alpha S$ random field, with $0 < \alpha < 2$, integral
representation $(\ref{eqn3})$, and
functions $\{f_t\}$ given by $(\ref{eqn4})$. Then, in the terminology
introduced in this section, we have the following: \\

\noindent 1. If $\{\phi_t\}_{t \in F}$ is not conservative then
\begin{equation}
\frac{1}{n^{p/\alpha}}M_n \Rightarrow c\,
Z_\alpha\label{M_n_rate_Gnoncons}
\end{equation}
for some $c\in (0,\infty)$, and $Z_\alpha$ as in
$(\ref{M_n_rate_noncons})$. In fact,
\[
c=\left(\frac{\mathcal{V}lC_\alpha}{2^p}\int_W (g^{\ast}(w))^\alpha
d\nu(w)\right)^{1/\alpha}\,,
\]
where $\mathcal{V}$ is given by $(\ref{limit_density})$, while
$g^\ast$ is given by $(\ref{defn_of_g^star})$ applied to the
dissipative part of the random field $\eqref{e:H.field}$, and
$C_\alpha$ is as in $(\ref{M_n_rate_noncons})$.\\

\noindent 2. If $\{\phi_t\}_{t \in F}$ is conservative then
\begin{equation}
\frac{1}{n^{p/\alpha}}M_n \probconv 0\,.\label{M_n_rate_Gcons}
\end{equation}
\end{thm}

\begin{proof} 1. Let $r_n$ be the left hand side of
$(\ref{suff_condn})$. Then we have,
\begin{eqnarray}
r_n  &\leq &  P\bigg(\mbox{for some }u \in H_n,
   \frac{|f_u(U_j^{(n)})|}{\max_{s \in H_n} |f_s(U_j^{(n)})|} >
   \epsilon,\; j=1,2 \bigg) \nonumber\\
   &\leq&
   |H_n|\left(\epsilon^{-\alpha}b_n^{-\alpha}\int_S|f(s)|^\alpha
   \mu(ds)\right)^2.    \label{bound_on_r_n}
\end{eqnarray}
The inequality $(\ref{bound_on_r_n})$ follows using the argument
given in Remark $4.2$ of \cite{samorodnitsky:2004a}. Since
$\{\phi_t\}_{t \in F}$ is not conservative, Proposition
$\ref{propn_b_n_Grate}$ yields that $b_n$ satisfies
$(\ref{b_n_rate_Gdiss})$. Hence by $(\ref{size_of_H_n})$ we get that
$(\ref{suff_condn})$ holds in this case. Since $b_n$ satisfies
$(\ref{b_n_rate_Gdiss})$ with $a$ given by $(\ref{form_of_a})$, we
get $(\ref{M_n_rate_Gnoncons})$ by Theorem $\ref{M_n_rate_thm}$.\\

\noindent 2. As in the proof of $(4.3)$ in
\cite{samorodnitsky:2004a} we can get a stationary $S\alpha S$
random field $\mathbf{Y}$ generated by a conservative
$\mathbb{Z}^d$-action such that $b_n^Y$ satisfies
$(\ref{inequality})$ as well as $(\ref{b_n_rate_Gcons})$ (this is
possible, for instance, by Example $\ref{ex_imp_for_thm}$ below).
Therefore, $(\ref{M_n_rate_Gcons})$ follows using the exact same
argument as in the proof of $(4.3)$ in \cite{samorodnitsky:2004a}.
\end{proof}

\begin{remark} \label{rk:p0}
\textnormal{The previous discussion asssumes that $p\geq 1$.
When $p=0$ (i.e. when $\mathbb{Z}^d/K$ is a finite group) the random
field takes only finitely many different values. Therefore,
the sequence $M_n$ remains constant after some stage and so
converges to the maximum of finitely many $X_t$'s, not an extreme
value random variable.}
\end{remark}


\section{Examples}\label{examples}

In this section we consider several examples of stationary $S\alpha S$
random
fields associated with conservative flows. As in the one-dimensional
case considered in \cite{samorodnitsky:2004a}, the idea is to exhibit
a variety of possible in this case behaviour.\\

\noindent The first example is parallel to examples 5.1 and 5.4 in
\cite{samorodnitsky:2004a}.


\begin{example} \label{ex_imp_for_thm} {\rm

Let the random field have an integral representation of the form
\begin{equation}
X_t \eqdef  \int_{\mathbb{R}^{\mathbb{Z}^d}} g_t\,
dM,\;\;\;t \in \mathbb{Z}^d \label{repn_of_subguassian}
\end{equation}
where  $M$ is a $S\alpha S$ random measure on
$\mathbb{R}^{\mathbb{Z}^d}$ whose
control measure $\mu$ is a probability measure under which the
projections $(g_t,\;t \in \mathbb{Z}^d)$ are i.i.d.
random variables, with a finite absolute $\alpha$th moment.

If $(g_t,\;t \in \mathbb{Z}^d)$ are i.i.d. standard normal random
variables under $\mu$, then, as in the one-dimensional case, one
sees that
\[
b_n^\alpha \sim (2d\log{n})^{\alpha/2}\,,
\]
the assumption $(\ref{inequality})$ in Theorem $\ref{M_n_rate_thm}$
fails, and $b_n^{-1}M_n$ converges to a nonextreme value limit. See
also Remark \ref{rk:p0} above.

On the other hand, if, under $\mu$,  $(g_t,\;t \in \mathbb{Z}^d)$ are
i.i.d. positive Pareto random variables with
\[
\mu(g_0>x)=x^{-\theta}\;\;\;\mbox{for }x \geq 1
\]
for some $\theta > \alpha$, then as in the one-dimensional case we see
that
\[
b_n  \sim c_{\alpha,\theta}^{1/\alpha}\,n^{d/\theta}\;\;\;\mbox{as
}n \rightarrow \infty\,,
\]
for some finite positive constant $c_{p,\theta}$, Theorem
\ref{M_n_rate_thm} applies, and $n^{-d/\theta}M_n$ converges to an
extreme value distribution and hence this example also shows that
the rate of growth of $M_n$ can be $n^\gamma$ for any $\gamma \in
(0, d/\alpha)$. Note that existence of such a process was needed in
the proof of $(\ref{M_n_rate_Gcons})$ in Theorem
$\ref{thm_M_n_Grate}$.}
\end{example}

Next is an example of an application of Theorem $\ref{thm_M_n_Grate}$.

\begin{example}
{\rm Suppose $d=3$, and define the $\mathbb{Z}^3$-action
$\{\phi_{(i,j,k)}\}$ on $S=\mathbb{R}\times \{-1,1\}$ as
\[
\phi_{(i,j,k)}(x,y)=(x+i+2j,(-1)^k y)\,.
\]
An action-invariant measure $\mu$ on $S$ is defined as the product
of the Lebesgue measure on $\mathbb{R}$ and the counting measure on
$\{-1,1\}$.

Take any $f \in L^\alpha(S)$ and define a stationary $S\alpha S$
random field
$\{X_{(i,j,k)}\}$ as follows
\[
X_{(i,j,k)}=\int_{\mathbb{R}\times \{-1,1\}}
f\big(\phi_{(i,j,k)}(x,y)\big)\, dM(x,y)\,,
\]
where $M$ is a $S\alpha S$ random measure on $\mathbb{R}\times
\{-1,1\}$ with control measure $\mu$. Note that the above
representation of $\{X_{(i,j,k)}\}$ is of the form $(\ref{eqn4})$
generated by a measure preserving conservative action with
$c_{(i,j,k)} \equiv 1$.}

{\rm In the notation of Section $\ref{sec_group_theory}$ we
have
\[
K=\{(i,j,k)\in \mathbb{Z}^3:\,i+2j=0 \mbox{ and }k \mbox{ is
even}\}\,,
\]
and so
\[
A \simeq \mathbb{Z}^3 / K \simeq
\mathbb{Z}\times\mathbb{Z}/2\mathbb{Z}\,,
\]
and
\[
F=\{(i,0,0): i\in \mathbb{Z}\}\,.
\]
In particular $p=1$ and $\{\phi_t\}_{t\in F}$ is dissipative. Hence
Theorem $\ref{thm_M_n_Grate}$ applies and says that
$\frac{1}{n^{1/\alpha}}M_n$
converges to an extreme value distribution.}
\end{example}

In all the examples we have seen so far, the action has a
conservative direction i.e there is $u \in \mathbb{Z}^d-\{0\}$ such
that $\{\phi_{nu}\}_{n \in \mathbb{Z}}$ is a conservative
$\mathbb{Z}$-action. The following example of a
$\mathbb{Z}^2$-action, suggested to us by
M.G. Nadkarni, lacks such a conservative direction. In a sense, this
example is ``less one-dimensional'' than the previous examples.
\begin{example}
\textnormal{Suppose that $d=2$, and define the action
$\{\phi_{(i,j)}\}_{i,j \in \mathbb{Z}}$ of $\mathbb{Z}^2$ on
$S=\mathbb{R}$ with $\mu=Leb$ by
\[
\phi_{(i,j)}(x)=x+i+j\sqrt{2}, \;\;\;\forall\,x \in \mathbb{R}\,.
\]
Clearly, this action is measure preserving and it does not have any
conservative direction. It is, however, well known that this action does not admit a wandering
set of positive Lebesgue measure, and hence is conservative. In fact, if we take the kernel
$f=I_{[0,1]}$ and define $\{X_{(i,j)}\}$ by $(\ref{eqn3})$ and
$(\ref{eqn4})$ with, say, $c_{(i,j)}\equiv 1$, then we have, for
all $n \geq 2$,
\[
b_n^\alpha =\mu\bigg( \bigcup_{0 \leq i,j \leq (n-1)}
\phi_{(i,j)}\big([0,1]\big)\bigg)=\mu\big([0,1+(n-1)(1+\sqrt{2})]\big)\,.
\]
So, $b_n\,\sim\,(1+\sqrt{2})^{1/\alpha}n^{1/\alpha}$ and, a simple
calculation shows that left hand side of $(\ref{suff_condn})$ is
bounded from above by
$b_n^{-2\alpha} (\mu \otimes \mu)(B_n)$ where
\[
B_n=\big\{(x,y)\in \mathbb{R}^2:\,-(n-1)(1+\sqrt{2}) \leq x,y \leq
1,\,|x-y| \leq 1\big\}\,.
\]
Since $(\mu \otimes \mu)(B_n)=O(n)$, $(\ref{suff_condn})$ holds and
hence
\[
\frac{1}{n^{1/\alpha}}M_n \Rightarrow
\big((1+\sqrt{2})C_\alpha\big)^{1/\alpha}Z_\alpha\,.
\]}

\end{example}

\noindent \textbf{Acknowledgement. }The authors are thankful to
Mahendra Ganpatrao Nadkarni and Laurent Saloff-Coste for a number of
useful discussions and to the anonymous referees for their comments.

\end{document}